\documentclass[final]{siamltex}
\usepackage[papersize={200mm,280mm},top=2cm, bottom=3cm, left=2cm , right=2cm]{geometry}
\usepackage{amsmath}
\usepackage{graphics}
\usepackage{latexsym}
\usepackage{amsfonts}
\usepackage{stmaryrd}
\usepackage{color}

\title{The Finite Horizon impulse control  Problem with arbitrary cost \newline functions : the Viscosity Solution Approach.}
\author{Brahim EL ASRI\footnotemark[1]\thanks{Universit\'e Ibn Zohr, Equipe. Aide \`a la decision, ENSA, B.P.  1136, Agadir, Maroc. e-mail: b.elasri@uiz.ac.ma}\
	\and  Sehail MAZID\footnotemark[2]\thanks{Universit\'e Ibn Zohr, Equipe. Aide \`a la decision,
		ENSA, B.P.  1136, Agadir, Maroc. e-mail: sehail.mazid@edu.uiz.ac.ma.}\
	}

\begin{document}
\maketitle
\newtheorem{theo}{Theorem}
\newtheorem{problem}{Problem}
\newtheorem{pro}{Proposition}
\newtheorem{cor}{Corollary}
\newtheorem{axiom}{Definition}
\newtheorem{rem}{Remark}
\newtheorem{lem}{Lemma}
\newcommand{\brm}{\begin{rem}}
\newcommand{\erm}{\end{rem}}
\newcommand{\beth}{\begin{theo}}
\newcommand{\eeth}{\end{theo}}
\newcommand{\bl}{\begin{lem}}
\newcommand{\el}{\end{lem}}
\newcommand{\bp}{\begin{pro}}
\newcommand{\ep}{\end{pro}}
\newcommand{\bcor}{\begin{cor}}
\newcommand{\ecor}{\end{cor}}
\newcommand{\be}{\begin{equation}}
\newcommand{\ee}{\end{equation}}
\newcommand{\beq}{\begin{eqnarray*}}
\newcommand{\eeq}{\end{eqnarray*}}
\newcommand{\beqa}{\begin{eqnarray}}
\newcommand{\eeqa}{\end{eqnarray}}
\newcommand{\dg}{\displaystyle \delta}
\newcommand{\cm}{\cal M}
\newcommand{\cF}{{\cal F}}
\newcommand{\cR}{{\cal R}}
\newcommand{\bF}{{\bf F}}
\newcommand{\tg}{\displaystyle \theta}
\newcommand{\w}{\displaystyle \omega}
\newcommand{\W}{\displaystyle \Omega}
\newcommand{\vp}{\displaystyle \varphi}
\newcommand{\ig}[2]{\displaystyle \int_{#1}^{#2}}
\newcommand{\integ}[2]{\displaystyle \int_{#1}^{#2}}
\newcommand{\produit}[2]{\displaystyle \prod_{#1}^{#2}}
\newcommand{\somme}[2]{\displaystyle \sum_{#1}^{#2}}
\newlength{\inter}
\setlength{\inter}{\baselineskip} \setlength{\baselineskip}{7mm}
\newcommand{\no}{\noindent}
\newcommand{\rw}{\rightarrow}
\def \ind{1\!\!1}
\def \R{I\!\!R}
\def \N{I\!\!N}
\def \cadlag {{c\`adl\`ag}~}
\def \esssup {\mbox{ess sup}}
\begin{abstract}
\normalsize{We consider stochastic impulse control problems when the impulses cost
functions  are  arbitrary. We use the dynamic programming principle and viscosity solutions approach to show that the value function is a unique viscosity solution for the associated Hamilton-Jacobi-Bellman equation (HJB) partial differential equation (PDE) of stochastic impulse control problems.}
\end{abstract}



\pagestyle{myheadings}
\thispagestyle{plain}

\section{Introduction}
\no

Impulse control problems form an important class of stochastic control problems. Finding a stochastic impulse control policy amounts to determining the sequence of random dates at which the policy is exercised and the sequence of impulses describing the magnitude of the applied policies, which maximizes a given reward function. It seems to be impossible to give an overview over all fields of application and all different variants that have been used. We only want to mention finance, e.g. cash management and portfolio optimization, see \cite{[K]} and \cite{[PS]},  control of an exchange rate by the Central Bank, see \cite{[MO]} and \cite{[CZ]}, and optimal forest management, see \cite{[W]}, \cite{[Al]} and the references therein.

In the literature, one finds several different approaches to tackle stochastic impulse control problems. One approach is to focus on solving the value function for the associated (quasi-)variational inequalities or
Hamilton-Jacobi-Bellman (HJB for short) integro-differential equations, and then establishing
the optimality of the solution by Verification Theorem (See \O ksendal and Sulem \cite{[OS]}).
Another approach is to characterize the value function of the control problem as a (unique) viscosity solution to the associated PDEs (See Lenhart \cite{[Le]}, Tang and Yong \cite{[TY]}). In dimension one other approaches, based on excessive mappings and iterated optimal stopping schemes (see \cite{[Al1],[AL],[Em],[HSZ]}).

The main objective of this paper is to study the problem of existence and uniqueness of a solution in viscosity
sense $V$ of the following system of partial differential equations with obstacles which depend on the solution:
\begin{equation} \label{eq:HJB0}\left\{
\begin{array}{c}
min\big\{-\partial_t V(t,x)-\mathcal{L}V(t,x)-f(t,x),V(t,x)-\sup\limits_{\xi\in U}[V(t,x+\xi)-c(t,x,\xi)]\big\}=0
\qquad \\
V(T,x)=g(x).\qquad\qquad\qquad\qquad\qquad\qquad\qquad\qquad\qquad\qquad\qquad\qquad\qquad\qquad\qquad
\end{array}
\right.
\end{equation}
Where $\mathcal{L}$ is the second-order local operator
$$\mathcal{L}V= \langle b,\nabla _x V\rangle +\cfrac{1}{2}tr[\sigma\sigma^*\nabla^2_xV].$$

In a way, the system (\ref{eq:HJB0}) is Bellman system of equations associated with the control impulse with utility functions $f$, terminal payoff $g$ and impulse cost given by $c$.

Amongst the papers which consider the same problem as ours, and
in the framework of viscosity solutions approach, the most elaborated works are certainly the ones by
Tang and Yong \cite{[TY]}, on the one hand, and by Seydel \cite{[Src]}. In \cite{[Src]}, the author restricted controls to be only Markov controls. This Markovian assumption simplifies the proof of dynamic programming principle  significantly. In the case of non-Markov controls, \cite{[TY]} show existence and uniqueness of a solution for (\ref{eq:HJB0}). Nevertheless the paper suffers from two facts: (i) the costs of the impulses to be decreasing in time. (ii) the costs of the impulses do not depend on the state variable. For the first note of \cite{[TY]} we can easily adapt the methods in \cite{[AZ]} or \cite{[EM1]}  in order to avoid the monotonicity
condition. The second issue of \cite{[TY]}, i.e. considering the case when $c$ depending also on $x$, was right now, according to our knowledge, an open problem. However, in that latter case, we face two main difficulties. First, it is not clear how to show the regularity of the value function, together with the gain functional, therefore we can not prove the dynamic programming principle. The second one is related to the obtention of the comparison of sub- and super-solutions of system (\ref{eq:HJB0}) which plays an important role in our study.

Therefore the main objective of our work, and this is the
novelty of the paper, is to show existence and uniqueness of a solution in viscosity sense for the system when the function $c$ is continuous depending also on $x$.  To derive
these results, we first study $(V^n)_n$, where $V^n$  is the value function from $t$ to $T$, when the system  only at most $n$ interventions after $t$ are allowed, we give the  dynamic programming principle for $V^n$ and we show that it is a continuous viscosity solution to
 \begin{equation}\label{Vn}\left\{
\begin{array}{c}
\min\big\{-\cfrac{\partial V^n}{\partial t}-\mathcal{L}V^{n}-f,V^n-\sup\limits_{\xi\in U}[V^{n-1}(t,x+\xi)-c(t,x,\xi)]\big\}=0
\qquad [0,T)\times\mathbb{R}^n, \\
V^n(T,x)=g(x)\qquad\qquad\qquad\qquad\qquad\qquad\qquad\qquad\qquad\qquad \qquad\qquad \forall x\in\mathbb{R}^n.
\end{array}
\right.
\end{equation}
Further, we  obtain that the $u.s.c.$ envelope $V^*$ of the value function, is a viscosity sub-solution of  (\ref{eq:HJB0}), and the $l.s.c.$ envelope $V_*$ of the value function, is a viscosity super-solution of  (\ref{eq:HJB0}). Finally, by comparison principle we obtain that the value function of the impulse control problem is a unique continuous viscosity solution to (\ref{eq:HJB0}).

This paper is organized as follows: In Section 2, we formulate the problem and we give the related definitions. In Section 3, we shall introduce the impulse control problem and  we study the problem in which the controller only can intervene finitely many times. Further we prove the dynamic programming principle. Section 4, is devoted to the connection between impulse control problem and Hamilton-Jacobi-Bellman equation. In Section 5, we show that the solution of Hamilton-Jacobi-Bellman equation is unique in the subclass of bounded functions.\\
\section{Assumptions and formulation of the problem}
\no

Throughout this paper $T$ (resp. $n,\, d$) is a fixed real (resp.
integers) positive numbers. $|.|$ will denote the canonical Euclidian norm on $\mathbb{R}^n$, and  $\langle .|.\rangle$ the corresponding inner product. Let us assume the following assumptions:  \medskip

\no\textbf{[H1]} $b:[0,T]\times \mathbb{R}^{n}\rightarrow
\mathbb{R}^{n}$ and $\sigma :[0,T]\times \mathbb{R}^{n}\rightarrow
\mathbb{R}^{n\times d}$ be two bounded continuous functions for which there
exists a constant $C>0$ such that for any $t\in \lbrack 0,T]$ and
$x,x^{\prime }\in \mathbb{R}^{n}$
\begin{equation}
|\sigma (t,x)-\sigma (t,x^{\prime })|+|b(t,x)-b(t,x^{\prime })|\leq
C|x-x^{\prime }|.
 \label{eqs}
\end{equation}%
\textbf{[H2]}  $f:[0,T]\times\mathbb{R}^{n} \rightarrow \mathbb{R}$
is uniformly continuous and bounded on $[0,T]\times\mathbb{R}^{n}$. $g:\mathbb{R}^{n} \rightarrow\mathbb{R}$ is uniformly continuous and bounded on $\mathbb{R}^n$.\\
\textbf{[H3]} The cost function $c:[0,T]\times\mathbb{R}^n\times U \rightarrow \mathbb{R}$ is measurable and uniformly continuous. Furthermore
\begin{equation}\label{cout} \underset{[0,T]\times\mathbb{R}^n\times U}{inf}c \geq k,
\end{equation}
where $k>0$. Moreover,
\begin{equation} \label{cout1}
  c(t,x,\xi_1+\xi_2)\leq c(t,x,\xi_1)+c(t,x,\xi_2),\\
\end{equation}
for every $(t,x)\in [0,T]\times\mathbb{R}^n, \xi_1,\xi_2\in U.$
\\
\textbf{[H4]} For any $x\in\mathbb{R}^n$ and $\xi\in U$
\begin{equation}\label{Terminal}
\sup\limits_{\xi\in U}[g(x+\xi)-c(T,x,\xi)]\leq g(x).
\end{equation}
Where $U$ is a  compact subset of $\mathbb{R}^{n}$.
\begin{rem}
Assumption (\textbf{H3})  and (\textbf{H4}) simply provides the classical framework for the study of impulse problems. (\ref{cout1}) ensures that multiple impulses occurring at the same time are suboptimal. (\ref{Terminal}) ensures the nonoptimality of a impulse at maturity. $\Box$
\end{rem}

We now consider the HJB equation:
\begin{equation} \label{eq:HJB}\left\{
\begin{array}{c}
min\big\{-\cfrac{\partial V}{\partial t}(t,x)-\mathcal{L}V(t,x)-f(t,x),V(t,x)-\mathcal{M} V(t,x)\big\}=0
\qquad [0,T)\times\mathbb{R}^n, \\ \\
V(T,x)=g(x)\qquad\qquad\qquad\qquad\qquad\qquad\qquad\qquad\qquad\qquad\qquad\qquad \forall x\in\mathbb{R}^n.
\end{array}
\right.
\end{equation}
Or equivalently,
\begin{equation}\left\{
\begin{array}{c}
-\displaystyle\cfrac{\partial V}{\partial t}(t,x)-\mathcal{L}V(t,x)-f(t,x)\geq 0\qquad\qquad\qquad\qquad\qquad\qquad\;[0,T)\times\mathbb{R}^n \\
V(t,x)\geq\mathcal{M} V(t,x)\qquad\qquad\qquad\qquad\qquad\qquad\qquad\qquad\qquad\quad[0,T)\times\mathbb{R}^n\\
\big(-\displaystyle\cfrac{\partial V}{\partial t}(t,x)-\mathcal{L}V(t,x)-f(t,x)\big)\big(V(t,x)-\mathcal{M} V(t,x)\big)=0\qquad[0,T)\times\mathbb{R}^n\\  V(T,x)=g(x)\qquad\qquad\qquad\qquad\qquad\qquad\qquad\qquad\qquad\qquad\qquad\forall x\in\mathbb{R}^n.
\end{array}
\right.
\end{equation}
Where $\mathcal{L}$ is the second-order local operator
   $$\mathcal{L}V= \langle b,\nabla _x V\rangle +\cfrac{1}{2}tr[\sigma\sigma^*\nabla^2_xV],$$
and the nonlocal operator $\mathcal{M}$  is given by
$$\mathcal{M}V(t,x)=\sup\limits_{\xi\in U}[V(t,x+\xi)-c(t,x,\xi)],$$
for every $(t,x)\in[0,T]\times\mathbb{R}^n.$

The main objective of this paper is to focus on the existence and uniqueness of
the solution in viscosity sense of (\ref{eq:HJB}) whose definition
is:
\begin{definition}
A lower (resp., upper) semicontinuous function $V:[0,T]\times\mathbb{R}^n\to\mathbb{R}$ is called a viscosity supersolution (resp., subsolution) to the HJB equation
(\ref{eq:HJB}) if
\begin{itemize}
  \item for every $(t_0,x_0) \in[0,T)\times\mathbb{R}^n$ and  any function $\phi \in C^{1,2}([0,T)\times\mathbb{R}^n)$, such that $(t_0,x_0)$ is a local minimum (resp., maximum)  of $V-\phi$, we have:
\begin{equation}
\begin{array}{ll}\label{visc1}
 \min\Big[-\cfrac{\partial \phi}{\partial t}(t_0,x_0)-\mathcal{L}\phi(t_0,x_0)-f(t_0,x_0),\\ \qquad
 \qquad V(t_0,x_0)-\mathcal{M} V(t_0,x_0)\Big]\geq 0\qquad (\text{resp.,}\leq 0).
\end{array}
\end{equation}
  \item for every $x\in \mathbb{R}^n$ we have
\begin{equation}
\begin{array}{ll}\label{visc2}
V(T,x)\geq g(x) \qquad (\text{resp.,}\leq 0).
\end{array}
\end{equation}
\end{itemize}

A bounded function $V:[0,T)\times\mathbb{R}^n\to\mathbb{R}$  is a viscosity solution to the HJB
equation (\ref{eq:HJB}) if its lower semicontinuous envelope $V_*$ is a viscosity supersolution and
its upper semicontinuous envelope $V^*$ is a viscosity subsolution. $V_*$ and $V^*$ are given by
$$V_*(t,x):=\liminf\limits_{(s,y)\to(t,x),s<T} V(s,y),\qquad\qquad V^*(t,x):=\limsup\limits_{(s,y)\to(t,x),s<T}V(s,y)$$
for every $(t,x)\in[0,T]\times\mathbb{R}^n$. $\Box$
\end{definition}

There is an equivalent formulation of this definition (see e.g. \cite{[TY]}) which we give because it will be
useful later. So firstly, we define the notions of superjet and subjet of a function $V$.
\begin{definition}
Let $V: [0,T]\times\mathbb{R}^n,$ be a lower (resp. upper) semicontinuous function, $(t,x)$ an element of $(0,T)\times\mathbb{R}^n$ and finally $S_n$ the set of $n\times n$ symmetric matrices. We denote by $J^{2,+}V(t,x)$ (resp. $J^{2,-}V(t,x))$, the superjets (resp. the subjets) of $\,V$
at $(t,x)$, the set of triples $(p,q,X)\in\mathbb{R}\times\mathbb{R}^n\times S_n$ such that:
\begin{equation*}
\begin{array}{lll}
V(s,y)\leq V(t,x)+p(s-t)+\langle q,y-x\rangle+\displaystyle\frac{1}{2}\langle X(y-x),y-x\rangle+o(|s-t|+|y-x|^2)
\end{array}
\end{equation*}
\begin{equation*}
\begin{array}{lll}
\big(\text{resp.}\\ V(s,y)\geq  V(t,x)+p(s-t)+\langle q,y-x\rangle+\displaystyle\frac{1}{2}\langle X(y-x),y-x\rangle+o(|s-t|+|y-x|^2)\big).\Box
\end{array}
\end{equation*}
\end{definition}
Note that if $\phi-V$ has a local maximum (resp. minimum) at $(t, x)$, then we obviously have:
$$(D_t\phi(t,x),D_x\phi(t,x),D^2_{xx}\phi(t,x))\in J^{2,-}V(t,x)\;(\text{resp}.\; J^{2,+}V(t,x)).\Box$$

We now give an equivalent definition of a viscosity solution of HJB equation(\ref{eq:HJB}):
\begin{definition}
Let $V$ be a lower (resp. upper)  continuous function defined on $[0,T]\times\mathbb{R}^n$. Then V is a viscosity supersolution (resp., subsolution) to the HJB equation
(\ref{eq:HJB}) if and only if
\begin{itemize}
\item for every $(t, x)\in [0,T)\times\mathbb{R}^n$ and $(p,q,X)\in J^{2,-}V(t,x)$\, (resp. $J^{2,+}V(t,x))$,
\begin{equation}
\begin{array}{ll}
min\Big[-p-\langle b(t,x),q\rangle-\displaystyle\frac{1}{2}Tr\big[(\sigma\sigma^*)(t,x)X\big] -f(t,x),V(t,x)-\mathcal{M}V(t,x)\Big]\geq 0 \qquad (\text{resp.}\leq0).
\end{array}
\end{equation}
 \item for every $x\in \mathbb{R}^n$ we have
\begin{equation}
\begin{array}{ll}\label{visc22}
V(T,x)\geq g(x) \qquad (\text{resp.}\leq 0).\qquad\Box
\end{array}
\end{equation}
\end{itemize}
\end{definition}

As pointed out previously we will show that system (\ref{eq:HJB})
has a unique solution in viscosity sense. This system is the
deterministic version of the stochastic impulse control problem
will describe briefly in the next section.
\section{The impulse control problem}
\subsection{Setting of the problem}
\subsubsection{Probabilistic setup}
\no

We work on a time horizon $[0,T]$, where $0 < T <\infty$. Let ($\Omega,\mathcal{F},\mathbb{P})$ is a fixed probability space on which is defined a
standard $d$-dimensional Brownian motion $W=(W_t)_{t\leq T}$, whose natural filtration is
$(\mathcal{F}^0_{t}:=\sigma\{W_s;s\leq t\})_{0\leq t\leq T}.$ We denote by $\mathbb{F} =(\mathcal{F}_{t})_{t\leq T}$ the completed filtration of $(\mathcal{F}^0_{t})_{t\leq T}$ with the $\mathbb{P}$-null sets of $\mathcal{F}$. The expectation operator with respect to $\mathbb{P}$ is denoted by $\mathbb{E}$, and the indicator function of a set or event $A$ is written as $\ind_{A}$. The notation $a.s.$ stands for almost-surely.
\subsubsection{Impulse control definitions}
\no

The following data for the impulse control problem are given:
\begin{enumerate}
	\item A spaces of control
	actions $U\subset \mathbb{R}^n$, where $U$ is a compact.
	\item A reward received at time $T$, which is modelled by an $\mathcal{F}_{T}$-measurable real-valued random \\ variable $g$.
	\item A running reward, which is represented by a real-valued
	adapted process $f=(f(t,x))_{(t,x)\in[0,T]\times\mathbb{R}^n}$.
	 \item A cost of the intervention $\xi\in U$, which is modelled by a real-valued adapted
	 process \\ $c=(c(t,x,\xi))_{(t,x)\in[0,T]\times\mathbb{R}^n}$.
\end{enumerate}

Define the concept of impulse control as follows:
\begin{definition}
An impulse control $u=\sum\limits_{m\geq1}\xi_{m}\ind_{[\tau_{m},T]}$  on $[t,T]\subset\mathbb{R}^{+}=[0,+\infty)$, is such that:
\begin{itemize}
	\item  $(\tau_{m})_{m}$, the action times, is a sequence of $\mathbb{F}$-stopping
	times, valued in $[t,T]\cup \{+\infty\}$ such that $\mathbb{P}$-a.s. $\tau_{m}\leq \tau_{m+1}$.
	\item $(\xi_{m})_{m}$, the actions, is a sequence of $U$-valued  random variables, where each $\xi_{m}$  is $\mathcal{F}_{\tau_{m}}$-measurable.
\end{itemize}
\end{definition}

We denote by $\mathcal{U}$ the set of processes $u(\cdot)$.\qquad$\Box$

Let $t\in[0,T]$ be the initial time  and $x \in \mathbb{R}^n$ the initial state. Then, given the impulse control $u$ on $[t,T]$, a stochastic process
$(X_s)_{s\geq 0}$ follows a stochastic differential equation,
\begin{equation}
\label{az}
\begin{array}{ll}
 X_{s}=x+\integ{t}{s}b(r,X_{r})dr+\integ{t}{s}\sigma(r,X_{r})dW_{r}
+\sum\limits_{m\geq 1}\xi_{m}\ind_{[\tau_{m},T]}(s)\qquad s\geq t.
\end{array}
\end{equation}
Note that the assumption (H1)  ensure, for any $(t,x)\in[0,T] \times \mathbb{R}^n$, the existence and uniqueness of a solution $X^{t,x,u}=\{X^{t,x,u}_s, t\leq s\leq T\}$ to the SDE (\ref{az}) (see \cite{[RY]} for more details).

Let $u=\sum_{m\geq1}\xi_{m} \ind_{[\tau_{m},T]}$  be an  impulse control  on $[t,T]$, and let $\tau \leq \sigma $ be two $[t,T]$-valued $\mathbb{F}$-stopping times. Then we define the restriction $u_{[\tau,\sigma]}$ of the impulse control $u$ by:
\begin{equation} u_{[\tau,\sigma]}(s)=\sum\limits_{m\geq1}\xi_{\mu_{t,\tau}(u)+m}\ind_{[\tau_{\mu_{t,\tau}(u)+m}\leq s\leq \sigma]}(s), \quad\quad\quad \tau \leq s \leq \sigma,
\end{equation}
$\mu_{t,\tau}$ is the number of impulses up to time $\tau$, i.e.,
$\mu_{t,\tau}(u):=\sum\limits_{m\geq 1}\ind_{[\tau_m\leq \tau]}.$

The stochastic control problem is to
\begin{equation}
\text{(Problem)}\quad\text{maximize}\quad J(t,x,u)\qquad \text{over all}\; u\in\mathcal{U}
\end{equation}
subject to (\ref{az}) with
\begin{equation}
\label{LAA}
J(t,x,u):=\mathbb{E}\bigg[\integ{t}{T}f(s,X_{s}^{t,x,u})ds-
\sum\limits_{m\geq 1} c(\tau_{m},X_{\tau_{m}},\xi_{m}) \ind_{[\tau_{m}\leq T]}+g(X_{T}^{t,x,u})\bigg].
\end{equation}
Here we denote V for the associated value function:
\begin{equation}\label{value function}
(\text{Value Function})\quad V(t,x)=\sup_{u\in\mathcal{U}}J(t,x,u).
\end{equation}

We first observe that the optimal impulse problem over $\mathcal{U}$ can be restricted
to the consideration of finite number of impulses $\mathcal{U}_{t,T}$, where $\mathcal{U}_{t,T}:=\{u=\sum_{m\geq1}\xi_{m} \ind_{[\tau_{m},T]}\in\mathcal{U}\;|\quad \mathbb{P}(\tau_k<T,\;\forall k\geq 1)=0\}. $
\begin{proposition}
Under (\textbf{H1}), (\textbf{H2}) and (\textbf{H3}) the supremum of $J$ over $\mathcal{U}$ coincides
with the one of $J$ over $\mathcal{U}_{t,T}$, that is
\begin{equation}
\sup_{u\in\mathcal{U}}J(t,x,u) = \sup_{u\in\mathcal{U}_{t,T}}J(t,x,u),\qquad\quad  (t,x)\in[0,T]\times\mathbb{R}^n.
\end{equation}
\end{proposition}
$Proof.$  If $u$ does not belong to $\mathcal{U}_{t,T}$, then $J(t,x,u) = -\infty$. Indeed, introduce $A:=\{\omega\in\Omega\,|\,\tau_n(\omega)< T,\;\forall n\geq 1\}$ and $A^c$ be its complement. Since $u\notin\mathcal{U}_{t,T}$, then $\mathbb{P}(A)> 0$. Using  the boundedness of $f$ and $g$, we deduce there exists
a constant $C> 0$ such that
\begin{equation}
J(t,x,u)\leq C-\mathbb{E}\big[\ind_{A}\big\{\sum_{m\geq 1} c(\tau_{m},X_{\tau_{m}},\xi_{m})\big\}+\ind_{A^c}\big\{\sum_{m\geq 1} c(\tau_{m},X_{\tau_{m}},\xi_{m})\big\}\big]=-\infty,
\end{equation}
since for any $(t,x)\in[0,T]\times\mathbb{R}^n$ and $\xi\in U, c(t,x,\xi)\geq k$, and directly deduce that $\sup_{u\in\mathcal{U}}J(t,x,u) = \sup_{u\in\mathcal{U}_{t,T}}J(t,x,u).$ \qquad $\Box$
\subsection {When the controller can intervene finitely many times}
\no

For $n=0,1,...$ let $U_{t,T}^n=\{u\in\mathcal{U}_{t,T}\quad \text{such that} \quad\tau_{n+1}=+\infty\}$. In other words, $U_{t,T}^n$ is the set of all  controls with at most $n$ interventions. Next let us define
\begin{equation}
V^n(t,x)=\sup_{u\in\mathcal{U}_{t,T}^n}J(t,x,u).
\end{equation}
\begin{lemma}
 The  sequence $(V^n)_{n\in\mathbb{N}}$ is increasing and converges on $[0,T]\times\mathbb{R}^n$ to the function $V$:
\begin{equation}
V^n(t,x)\nearrow V(t,x)\qquad\text{as}\quad n\nearrow +\infty
\end{equation}
for all $(t,x)\in [0,T]\times\mathbb{R}^n$.
\end{lemma}

$Proof$. Since $U_{t,T}^n\subset U_{t,T}^{n+1}\subset U_{t,T}$, it follows that $(V^n(t,x))_{n\in\mathbb{N}}$ is a non-decreasing sequence and
\begin{equation}
\forall (t,x)\in[0,T]\times\mathbb{R}^n\qquad \lim\limits_{n\to\infty} V^n(t,x)\leq V(t,x).
\end{equation}

Let us fix $t$ and $x$. For a given $\epsilon> 0$, there exists $u^\epsilon=\sum_{m\geq 1}\xi^\epsilon_{m}\ind_{[\tau^\epsilon_{m},T]}\in\mathcal{U}_{t,T}$ such that
\begin{equation}\label{epsiopt}
J(t,x,u^\epsilon)\geq V(t,x)-\epsilon.
\end{equation}
Define a control $u^n\in\mathcal{U}_{t,T}^n$ from  $u^\epsilon$ by setting $u^n=\sum_{m=1}^{n}\xi^\epsilon_{m}\ind_{[\tau^\epsilon_{m},T]}.$ Then
$X_s^{t,x,u^n}=X_s^{t,x,u^\epsilon}$ for all  $s\leq \tau_n$
and by positivity of the impulses costs,

\begin{eqnarray*}
J(t,x,u^\epsilon)-J(t,x,u^n) &\leq& \mathbb{E}\bigg[\bigg(\integ{t}{T}(f(s,X_{s}^{t,x,u^\epsilon})-f(s,X_{s}^{t,x,u^n}))ds+(g(X_{T}^{t,x,u^\epsilon})-g(X_{T}^{t,x,u^n}))\bigg)\ind_{\{\mu_{t,T}(u^\epsilon)>n\}}\bigg]\\ &\leq& C\mathbb{P}(\mu_{t,T}(u^\epsilon)>n)^{\frac{1}{2}},
\end{eqnarray*}
by Cauchy-Schwarz inequality and the boundedness of $f$ and $g$.  Hence letting $n \to \infty$,
and since $\mu_{t,T}(u^\epsilon) < \infty$ $a.s.$, we obtain
\begin{eqnarray*}
  \liminf\limits_{n\to\infty}J(t,x,u^n)&=&\liminf\limits_{n\to\infty}\mathbb{E}\bigg[\integ{t}{T}f(s,X_{s}^{t,x,u^n})ds-
\sum\limits_{m=1}^{n} c(\tau_{m},X_{\tau_{m}},\xi_{m}) \ind_{[\tau_{m}\leq T]}+g(X_{T}^{t,x,u^n})\bigg]\\
&\geq & J(t,x,u^\epsilon).
\end{eqnarray*}
Therefore, by (\ref{epsiopt})  we get
\begin{equation}
\liminf\limits_{n\to\infty}V^n(t,x)\geq \liminf\limits_{n\to\infty}J(t,x,u^n)\geq V(t,x)-\epsilon.
\end{equation}
Then by arbitrariness of $\epsilon$ we get  the required assertion.\qquad$\Box$

The following lemma is a direct consequence of the positivity of $c$ and the boundedness of $f$ and $g$. As such, its proof is omitted.
\begin{lemma}\label{boun}
Under the standing assumptions (\textbf{H1}), (\textbf{H2}) and (\textbf{H3}), the value function is bounded.
\end{lemma}
\begin{proposition}\label{cont} Under (\textbf{H1}), (\textbf{H2}) and (\textbf{H3}) we have
\begin{itemize}	
\item[(i)] $V^n$  is continuous in $[0,T]\times\mathbb{R}^n$
\item[(ii)] $V$ is  lower semicontinuous in $[0,T]\times\mathbb{R}^n$
\item[(iii)] $J(\cdot,u)$ is continuous in $[0,T]\times\mathbb{R}^n$ for all $u\in\mathcal{U}_{0,T}^n$.
\end{itemize}
\end{proposition}
$Proof.$
We only prove (i). The two other ones follow from similar arguments. Let $\epsilon>0$ , $t'\in]t,T],$ and $x'\in B(x,\epsilon),$ then there exist $u_\epsilon\in\mathcal{U}_{t,T}^n$ such that
\begin{equation}
V^n(t,x)-V^n(t',x')\leq J(t,x,u_\epsilon)-J(t',x',\bar{u}_\epsilon)+\epsilon,
\end{equation}
where $\bar{u}_\epsilon\in\mathcal{U}^n_{t',T}$ will be chosen later. In particular, suppose that $u_\epsilon:=\sum_{m=1}^{ n}\xi^\epsilon_{m}\ind_{[\tau^\epsilon_{m},T]}\in\mathcal{U}^n_{t,T}$, then define  $\bar{u}_\epsilon$ from $u_\epsilon$ by setting
$\bar{u}_\epsilon:=\sum_{m=1}^{ n}\xi^\epsilon_{m} \ind_{[t'\vee\tau^\epsilon_{m},T]}$.
Therefore we have
\begin{eqnarray*}
V^n(t,x)-V^n(t',x')&\leq& \mathbb{E}\big[\integ{0}{T}f(s,X_{s}^{t,x,u_\epsilon})\ind_{[s\geq t]}ds-
\sum\limits_{1\leq m\leq n} c(t\vee\tau^\epsilon_{m},X^{t,x}_{t\vee\tau^\epsilon_{m}},\xi^\epsilon_{m}) \ind_{[\tau^\epsilon_{m}\leq T]}+g(X_{T}^{t,x,u^\epsilon})\\ & &
-\integ{0}{T}f(s,X_{s}^{t',x',\bar{u}_\epsilon})\ind_{[s\geq t']}ds+
\sum\limits_{1\leq m\leq n} c(t'\vee\tau^\epsilon_{m},X^{t',x'}_{t'\vee\tau^\epsilon_{m}},\xi^\epsilon_{m}) \ind_{[\tau_{m}^\epsilon\leq T]}-g(X_{T}^{t',x',\bar{u}_\epsilon})\big]+\epsilon\\
&\leq& \mathbb{E}\big[\integ{0}{T}|f(s,X_{s}^{t,x,u_\epsilon})-f(s,X_{s}^{t',x',\bar{u}_\epsilon})|\ind_{[s\geq t']}ds+\integ{0}{T}|f(s,X_{s}^{t,x,u_\epsilon})|\ind_{[t\leq s < t']}\\ & & +n\max\limits_{1\leq m\leq n}\sup\limits_{s\leq T}|c(t'\vee s,X^{t',x'}_{t'\vee s},\xi^\epsilon_{m}) -c(t\vee s,X^{t,x}_{t\vee s},\xi^\epsilon_{m})|
+|g(X_{T}^{t,x,u^\epsilon})-g(X_{T}^{t',x',\bar{u}_\epsilon})|\big]+\epsilon.
\end{eqnarray*}
Moreover, by a standard estimate for the SDE applying It\^o's formula to
$|X_{s}^{t,x,u_\epsilon}-X_{s}^{t',x',\bar{u}_\epsilon}|^2$ and using Gronwall's lemma, we then obtain from the Lipschitz condition on $b$ and  $\sigma$
\begin{equation*}
\mathbb{E}|X_{s}^{t,x,u_\epsilon}-X_{s}^{t',x',\bar{u}_\epsilon}|\leq C(|x-x'|+|t-t'|^{\frac{1}{2}}) \qquad \forall s \in[t',T]
\end{equation*}
taking the limit as $(t,x)\to (t',x')$, and using the uniform continuity of $f$, $g$ and $c$ to obtain:
$$\limsup\limits_{(t,x)\to (t',x')}V^n(t,x)\leq V^n(t',x')+\epsilon.$$
As $\epsilon$ is arbitrary then sending $\epsilon\to 0$ to obtain:
 $$\limsup\limits_{(t,x)\to (t',x')}V^n(t,x)\leq V^n(t',x').$$
Therefore $V^n$ is upper semi-continuous. In a similar way we can prove that
$$\liminf\limits_{(t,x)\to (t',x')}V^n(t,x)\geq V^n(t',x').$$
Therefore $V^n$ is continuous.\qquad$\Box$
\subsection{Dynamic programming principle}
\no

The rigorous connection between $V$ and HJB equation passes through the dynamic programming principle (DPP). We begin with the following two lemmas.
\begin{lemma}
	Let $u=\sum_{m=1}^{n}\xi_{m} \ind_{[\tau_{m},T]}$ be a nearly optimal control of $V^n$, then there exist a positive
	constant $C$ which does not depend on $t$ and $x$ such that:
	\begin{equation}\label{coutpro}
	\mathbb{E}\bigg[\sum\limits_{m=1}^{n} c(\tau_{m},X_{\tau_{m}},\xi_{m}) \ind_{[\tau_{m}\leq T]}\bigg]\leq C.
	\end{equation}
	We denote by $\widehat{\mathcal{U}}^n_{t,T}$ the set of impulse controls which satisfies the condition (\ref{coutpro}).\\
\end{lemma}
$Proof.$ Let us choose a nearly optimal control $u=\sum_{m=1}^{n}\xi_{m} \ind_{[\tau_{m},T]}\in \mathcal{U}_{t,T}^n$ such that,
\begin{equation}
\mathbb{E}\bigg[\integ{t}{T}f(s,X_{s}^{t,x,u})ds-
\sum\limits_{m=1}^n c(\tau_{m},X_{\tau_{m}},\xi_{m}) \ind_{[\tau_{m}\leq T]}+g(X_{T}^{t,x,u})\bigg]\geq V^n(t,x)-1.
\end{equation}
Since $V^n, f$ and $g$ are bounded, then we obtain the desired result.\qquad$\Box$
\begin{lemma}
Under (\textbf{H1}), (\textbf{H2}) and (\textbf{H3})
we have
\begin{equation}
V^n(t,x)=\sup_{u\in\bar{\mathcal{U}}_{t,T}^n}J(t,x,u),
\end{equation}
for evry $(t,x)\in[0,T)\times\mathbb{R}^n$, where $\bar{\mathcal{U}}_{t,T}^n$ contain all the impulse controls
in $\widehat{\mathcal{U}}^n_{t,T}$ which have no impulses at time t.
\end{lemma}
\\$Proof.$
Let $u\in\widehat{\mathcal{U}}^n_{t,T}\setminus\bar{\mathcal{U}}_{t,T}^n,$ then for every $\epsilon> 0$, we have to prove that there exist $\bar{u}\in\bar{\mathcal{U}}_{t,T}^n$ such that
$$ |J(t,x, u) - J(t,x,\bar{u})|\leq \epsilon.$$
We only consider the case in which $u$  has  a single impulse at time $t$. (In the case of multiple impulses at time $t$, we using condition (\ref{cout1}), that we can
reduce this case with only a single impulse at time $t$).
Then there exist a $[t,T]-$valued $\mathbb{F}-$stopping times $\tau$, with
$\mathbb{P}(\tau=t) > 0$ such that  $$u=\xi\ind_{[\tau,T]}+\hat{u},$$
where $\hat{u}=\sum_{m=2}^{n}\xi_m\ind_{[\tau_m,T]}\in\bar{\mathcal{U}}_{t,T}^n$ and $\xi$ is an $\mathcal{F}_\tau-$ measurable
$U-$valued random variable.
For every integer $k$  (large enough) we pose $\tau_k=\big(\tau+\frac{1}{k}\big)\ind_{[\tau=t]}+\tau\ind_{[\tau>t]}$.
Next, define the impulse control
$$u_k=\xi\ind_{[\tau_k,T]}+\hat{u}\in\bar{\mathcal{U}}_{t,T}^n.$$
Then $\tau_k\to\tau$ as $k\to\infty$, $\mathbb{P}-a.s$. Moreover, for
all $s\in(t,T],\; X^{t,x,u_k}_s\to  X^{t,x,u}_s$ as $k\to\infty$, $\mathbb{P}-a.s.$
\\On the other hand, we have
\begin{eqnarray*}
J(t,x,u)-J(t,x,u_k)=\mathbb{E}\big[\integ{t}{T}(f(s,X_{s}^{t,x,u})-f(s,X_{s}^{t,x,u_k}))ds-c(\tau,X^{t,x}_{\tau},\xi) +c(\tau_k,X^{t,x}_{\tau_k},\xi)\\
+g(X_{T}^{t,x,u})-g(X_{T}^{t,x,u_k})\big].\qquad\qquad\qquad\qquad\qquad\qquad\qquad
\end{eqnarray*}
Therefore, from the dominated convergence theorem, we deduce the existence of an integer $K\geq 1$ such
that
$$|J(t,x,u)-J(t,x,u_k)| \leq \epsilon\qquad \forall k\geq K.$$
Thus we get the thesis.\qquad$\Box$
\begin{theorem}
Given $0\leq t\leq r\leq T$, $x\in\mathbb{R}^n$, and a $[t,r]-$valued $\mathbb{F}$ stopping times $\tau$, we have
\begin{eqnarray*}
V^n(t,x)&=& \sup_{u\in\widehat{\mathcal{U}}^n_{t,T}} \mathbb{E}\big[\int_{t}^{\tau\wedge\tau_1}f(s,X_{s}^{t,x,u})ds+(V^{n-1}(\tau_1,X_{\tau_1})-
c(\tau_{1},X_{\tau_{1}},\xi_{1}) )\ind_{[\tau_{1}\leq \tau]}+ V^{n}(\tau,X_{\tau})\ind_{[\tau<\tau_{1}]}\big].
\end{eqnarray*}
\end{theorem}
$Proof.$ Let $\epsilon>0$, then there exists $u^{n,\epsilon}\in\bar{\mathcal{U}}_{s,T}^n$ and $u^{n-1,\epsilon}\in\bar{\mathcal{U}}_{s,T}^{n-1}$ such that
\be \left\{
\begin{array}{ll}\label{progeq1}
V^{n}(s,y)\leq J(s,y,u^{n,\epsilon})+\epsilon \\
V^{n-1}(s,y)\leq J(s,y,u^{n-1,\epsilon})+\epsilon.
\end{array}
\right. \ee
Now, for fixed $(s,y)\in[0,T]\times\mathbb{R}^n$, the continuity of $V^n$ and $J$  established in Proposition \ref{cont} imply that there exists $r^{(\epsilon,s,y)}$ such that
\be \left\{
\begin{array}{ll}\label{progeq2}
V^l(t',y')\leq V^l(s,y)+\epsilon\\
J^l(t',y',u^{l,\epsilon})\geq J(s,y,u^{l,\epsilon})-\epsilon,
\end{array}
\right.\qquad\qquad l=n,n-1 \ee
for all $(t',y')\in B(s,y;r^{(\epsilon,s,y)})$, where for $r>0$ and $(s,y)\in[0,T]\times\mathbb{R}^n$
\begin{equation}
B(s,y;r)=\{(t',y')\in[0,T]\times\mathbb{R}^n:t'\in(s-r,s],|y'-y|< r\}.
\end{equation}
Therefore, the family $\{B(s,y;r): (s,y)\in[0,T]\times\mathbb{R}^n, 0<r\leq r^{(\epsilon,s,y)}\}$ forms
an open covering of $(0,T]\times\mathbb{R}^n$. By the $Lindel\ddot{o}f$ covering Theorem (\cite{[Dug]}, Theorem 6.3 Chapter
VIII), there exists a countable sequence $(t_i,y_i,r_i)_{i\geq 1}$ in $[0,T]\times\mathbb{R}^n\times\mathbb{R}$ such that $\{B(t_i,y_i;r_i)\}_{i\geq 1}$ is a countable subcover of $(0,T]\times\mathbb{R}^n$. We set $u_{n,i}^\epsilon:=u^{n,(t_i,y_i),\epsilon}$ and $u^\epsilon_{n-1,i}:=u^{n-1,(t_i,y_i),\epsilon}$.\\ Now set $A_0:=\{T\}\times\mathbb{R}^n$, $C_{-1}:=\emptyset$ and define for all $i\in\mathbb{N}\cup \{0\}$
$$A_{i+1}:=B(t_{i+1},y_{i+1};r_{i+1})\setminus C_i,\qquad \text{where}\; C_i:=C_{i-1}\cup A_i, \qquad i\geq 0.$$
Under this construction, we have
\begin{equation}
(\tau\wedge\tau_1,X_{\tau\wedge\tau_1})\in\cup_{i\in\mathbb{N}\cup\{0\}}A_i \;\mathbb{P}-a.s.,\quad\text{and} \;\, A_i\cap A_j=\emptyset \quad \text{for}\quad i\ne j.
\end{equation}
Next the inequalities $(\ref{progeq1})$ and $(\ref{progeq2})$ yield that
\begin{equation}\label{progeq3}
V^l(t',y')\leq J(t',y',u^\epsilon_{l,i}) +3\epsilon \quad\text{for all}\; (t',y')\in A^i\;\,\text{and}\;\,  l=n,n-1.
\end{equation}
For any $k\in\mathbb{N}$, set $A^k:=\bigcup_{0\leq i\leq k}A_i$. Given $u\in\widehat{\mathcal{U}}_{t,T}^n$, we define
$$ u^k:=u\ind_{[t,\tau\wedge\tau_1]}+\ind_{(\tau\wedge\tau_1,T]}\bigg(u\ind_{(A^k)^c}+\sum\limits_{1\leq i\leq k}u^\epsilon_{n-1,i}\ind_{A_i}\ind_{[\tau_{1}\leq \tau]}+\sum\limits_{1\leq i\leq k}u^\epsilon_{n,i}\ind_{A_i}\ind_{[\tau_{1}> \tau]}\bigg).$$
Then
\begin{equation}
\begin{array}{ll}
J(t,x,u^k)=  \mathbb{E}\bigg[\integ{t}{\tau\wedge\tau_1}f(r,X_{r}^{t,x,u})dr-
c(\tau_{1},X_{\tau_{1}},\xi_{1})\ind_{[\tau_{1}\leq \tau]} +J(t,x,u)\ind_{(\Gamma(k))^c}\\ \\ +\sum\limits_{1\leq i,j\leq k}(J(\tau_1,X_{\tau_1}^{t,x,u},u^\epsilon_{n-1,i})\ind_{[\tau_{1}\leq \tau]}+J(\tau,X_{\tau}^{t,x,u},u^\epsilon_{n,i})\ind_{[\tau_{1}> \tau]})\ind_{\Gamma_j^i}\bigg].
\end{array}
\end{equation}
We deduce via $(\ref{progeq3})$ that
\begin{equation}\begin{array}{ll}\label{progeq4}
\mathbb{E}\bigg[\sum\limits_{1\leq i,j\leq k}(J(\tau_1,X_{\tau_1}^{t,x,u},u^\epsilon_{n-1,i})\ind_{[\tau_{1}\leq \tau]}+J(\tau,X_{\tau}^{t,x,u},u^\epsilon_{n,i})\ind_{[\tau_{1}> \tau]})\ind_{\Gamma_j^i}\bigg]\\ \\ \qquad\qquad\qquad\qquad\qquad\qquad\geq\mathbb{E}\big[(V^{n-1}(\tau_1,X_{\tau_1}^{t,x,u})\ind_{[\tau_{1}\leq \tau]}+V^{n}(\tau,X_{\tau}^{t,x,u})\ind_{[\tau_{1}> \tau]})\ind_{\Gamma(k)}\big]-3\epsilon,
\end{array}
\end{equation}
for every $k\geq1$. Letting $k\rightarrow\infty$, therefore,
$$\mathbb{E}\big[J(t,x,u)\ind_{(\Gamma(k))^c}\big]\longrightarrow 0$$
by dominated convergence and since $J(t,x,u)$ is bounded.  Moreover, monotone convergence yields
$$\mathbb{E}\big[(V^{n-1}(\tau_1,X_{\tau_1}^{t,x,u})\ind_{[\tau_{1}\leq \tau]}+V^{n}(\tau,X_{\tau}^{t,x,u})\ind_{[\tau_{1}> \tau]})\ind_{\Gamma(k)}\big]\longrightarrow \mathbb{E}\big[(V^{n-1}(\tau_1,X_{\tau_1}^{t,x,u})\ind_{[\tau_{1}\leq \tau]}+V^{n}(\tau,X_{\tau}^{t,x,u})\ind_{[\tau_{1}> \tau]})\big].$$
Therefore, we deduce the existence of an integer $k_0 \geq 1$ such that
\begin{equation*}
\begin{array}{ll}\label{progeq41}
J(t,x,u^{k_0})\geq  \mathbb{E}\bigg[\integ{t}{\tau\wedge\tau_1}f(s,X_{s}^{t,x,u})ds+(V^{n-1}(\tau_1,X_{\tau_1})-
c(\tau_{1},X_{\tau_{1}},\xi_{1}) )\ind_{[\tau_{1}\leq \tau]}+ V^{n}(\tau,X_{\tau})\ind_{[\tau<\tau_{1}]}\bigg]-4\epsilon.
\end{array}
\end{equation*}
The arbitrariness of $u$  and $\epsilon$ implies that
\begin{equation}
\begin{array}{ll}\label{progeq5}
V^n(t,x)\geq \sup\limits_{u\in\widehat{\mathcal{U}}^n_{t,T}}\mathbb{E}\bigg[\integ{t}{\tau\wedge\tau_1}f(s,X_{s}^{t,x,u})ds+(V^{n-1}(\tau_1,X_{\tau_1})-
c(\tau_{1},X_{\tau_{1}},\xi_{1}) )\ind_{[\tau_{1}\leq \tau]}+ V^{n}(\tau,X_{\tau})\ind_{[\tau<\tau_{1}]}\bigg].
\end{array}
\end{equation}
On the other hand for any $u\in\widehat{\mathcal{U}}^n_{t,T}$, we have
\begin{eqnarray*}
J(t,x,u)
&\leq& \mathbb{E}\bigg[\integ{t}{\tau\wedge\tau_1}f(s,X_{s}^{t,x,u})ds+(J(\tau_1,X_{\tau_1},u_{[\tau_1,T]})-
c(\tau_{1},X_{\tau_{1}},\xi_{1}) )\ind_{[\tau_{1}\leq \tau]}+ J(\tau,X_{\tau},u_{[\tau,T]})\ind_{[\tau<\tau_{1}]}\bigg]\\
&\leq& \mathbb{E}\bigg[\integ{t}{\tau\wedge\tau_1}f(s,X_{s}^{t,x,u})ds+(V^{n-1}(\tau_1,X_{\tau_1})-
c(\tau_{1},X_{\tau_{1}},\xi_{1}) )\ind_{[\tau_{1}\leq \tau]}+ V^{n}(\tau,X_{\tau})\ind_{[\tau<\tau_{1}]}\bigg].
\end{eqnarray*}
Taking supremum on both sides, we get
\begin{equation}
\begin{array}{ll}\label{progeq6}
V^n(t,x)\leq \sup\limits_{u\in\widehat{\mathcal{U}}^n_{t,T}}\mathbb{E}\bigg[\integ{t}{\tau\wedge\tau_1}f(s,X_{s}^{t,x,u})ds+(V^{n-1}(\tau_1,X_{\tau_1})-
c(\tau_{1},X_{\tau_{1}},\xi_{1}) )\ind_{[\tau_{1}\leq \tau]}+ V^{n}(\tau,X_{\tau})\ind_{[\tau<\tau_{1}]}\bigg].
\end{array}
\end{equation}
Then from (\ref{progeq5}) and (\ref{progeq6}) we deduce the thesis.\qquad$\Box$
\section{Hamilton-Jacobi-Bellman equation}
\no

In the present section we study HJB equation by means of viscosity solutions.
\begin{theorem}
$V^n$ is a continuous viscosity solution to
 \begin{equation}\label{eq:HJBn}\left\{
\begin{array}{c}
\min\big\{-\cfrac{\partial V^n}{\partial t}-\mathcal{L}V^{n}-f,V^n-\mathcal{M} V^{n-1}\big\}=0
\qquad [0,T)\times\mathbb{R}^n, \\ \\
V^n(T,x)=g(x)\qquad\qquad\qquad\qquad\qquad\qquad\qquad \forall x\in\mathbb{R}^n.
\end{array}
\right.
\end{equation}
\end{theorem}
$Proof.$ We know from Proposition \ref{cont} that $V^n$ is continuous on $[0,T]\times \mathbb{R}^n$, therefore
$V^n$ is equal to its lower semicontinuous envelope and to its upper semicontinuous envelope on $[0,T]\times \mathbb{R}^n$. We begin by proving that $V^n$ is a viscosity subsolution to (\ref{eq:HJBn}). Suppose $V^n-\phi$ achieves a local maximum in $[t_0,t_0+\delta)\times B(x_0,\delta)$ with $V^n(t_0,x_0)=\phi(t_0,x_0)$.
When $\tau = t$ in the dynamic programming principle for $V^n$, we have,
$V^n(t_0,x_0)\geq \mathcal{M} V^{n-1}(t_0,x_0)$. Then if $ V^n(t_0,x_0)= \mathcal{M} V^{n-1}(t_0,x_0)$, we already have the desired
inequality. Now suppose
$$V^n(t_0,x_0)-\mathcal{M} V^{n-1}(t_0,x_0)>2\epsilon>0,$$
we prove by contradiction that
\begin{equation}
\begin{array}{ll} -\cfrac{\partial \phi}{\partial t}(t_0,x_0)-\mathcal{L}\phi (t_0,x_0)-f(t_0,x_0) \leq 0.
\end{array}
\end{equation}
Suppose otherwise, i.e.,$\, -\cfrac{\partial \phi}{\partial t}(t_0,x_0)-\mathcal{L}\phi (t_0,x_0)-f(t_0,x_0)> 0$. Then without
loss of generality we can assume that $ -\cfrac{\partial \phi}{\partial t}(t,x)-\mathcal{L}\phi(t,x)-f(t,x) > 0$ and $V^n(t,x)-\mathcal{M} V^{n-1}(t,x)>\epsilon$ on $[t_0,t_0+\delta)\times B(x_0, \delta)$.\\
Define the
stopping time $\tau$ by
$$\tau=\inf\{t\in[t_0,T]:(t,X_t)\notin [t_0,t_0+\delta)\times B(x_0, \delta) \}\land T.$$
By It\^o's formula
\begin{equation}
\begin{array}{llll}
\label{tt}
\mathbb{E}[\phi(\tau,X_\tau^{t_0,x_0,u_0})]-\phi(t_0,x_0)=\mathbb{E}\bigg[\integ{t_0}{\tau}\bigg(\cfrac{\partial \phi}{\partial t}+\mathcal{L}\phi \bigg)(r,X_{r}^{t_0,x_0,u_0})dr\bigg],
\end{array}
\end{equation}
where $u_0$  is the control with no impulses.
\\Let $\epsilon_1>0$, using the dynamic programming principle between time $t_0$ and $\tau\wedge\tau_1$, we deduce the existence of a control $u^{\epsilon_1}\in\widehat{\mathcal{U}}^n_{t_0,T}$ such that
\begin{eqnarray*}
V^{n}(t_0,x_0)&\leq& \mathbb{E}\big[\int_{t_0}^{\tau\wedge\tau^\epsilon_1}f(r,X_{r}^{t_0,x_0,u^{\epsilon_1}})dr+(V^{n-1}(\tau^{\epsilon_1}_1,X^{t_0,x_0,u^{\epsilon_1}}_{\tau^{\epsilon_1}_1})-
c(\tau^{\epsilon_1}_{1},X^{t_0,x_0,u^{\epsilon_1}}_{\tau^{\epsilon_1}_{1}},\xi_{1}) )\ind_{[\tau^{\epsilon_1}_{1}\leq \tau]}\\ & &+ V^{n}(\tau,X^{t_0,x_0,u^{\epsilon_1}}_{\tau})\ind_{[\tau<\tau^{\epsilon_1}_{1}]}\big]+\epsilon_1
\\ &\leq& \mathbb{E}\big[\int_{t_0}^{\tau\wedge\tau^{\epsilon_1}_1}f(r,X_{r}^{t_0,x_0,u_0})dr+\mathcal{M} V^{n-1}(\tau^{\epsilon_1}_1,X^{t_0,x_0,u_0}_{\tau^\epsilon_1})\ind_{[\tau^{\epsilon_1}_{1}\leq \tau]}+ V^{n}(\tau,X^{t_0,x_0,u_0}_{\tau})\ind_{[\tau<\tau^{\epsilon_1}_{1}]}\big]+\epsilon_1
\\ &\leq& \mathbb{E}\bigg[\integ{t_0}{\tau\wedge\tau^{\epsilon_1}_1}f(r,X_{r}^{t_0,x_0,u_0})dr+V^{n}(\tau^{\epsilon_1}_1\wedge\tau,X_{\tau^{\epsilon_1}_1\wedge\tau}^{t_0,x_0,u_0})\bigg]-\epsilon.\mathbb{P}(\tau^{\epsilon_1}_1\leq\tau)+\epsilon_1.
\end{eqnarray*}
Therefore, without loss of generality, we only need to consider $u^{\epsilon_1}\in\widehat{\mathcal{U}}^n_{t_0,T}$ such that $\tau^\epsilon_1 > \tau$. Then
\begin{eqnarray*}
\phi(t_0,x_0)=V^{n}(t_0,x_0)&\leq& \mathbb{E}\bigg[\integ{t_0}{\tau}f(r,X_{r}^{t_0,x_0,u_0})dr+ V^{n}(\tau,X_{\tau}^{t_0,x_0,u_0})\bigg]+\epsilon_1
 \\ &\leq& \mathbb{E}\bigg[\integ{t_0}{\tau}f(r,X_{r}^{t_0,x_0,u_0})dr+ \phi(\tau,X_{\tau}^{t_0,x_0,u_0})\bigg]+\epsilon_1.
\end{eqnarray*}
Then from (\ref{tt}) and sending $\epsilon_1\rightarrow 0 $ we get
$$
 0\geq \mathbb{E}\bigg[\integ{t_0}{\tau}-\bigg(\cfrac{\partial \phi}{\partial t}+\mathcal{L}\phi+f\bigg)(r,X_{r}^{t_0,x_0,u_0})dr\bigg],
$$
which is a contradiction. Therefore,
$$\min\big\{-\cfrac{\partial \phi}{\partial t}-\mathcal{L}\phi-f,V^n-\mathcal{M} V^{n-1}\big\}(t_0,x_0)\leq0,$$
which is the subsolution property. The supersolution  property is proved analogously.\qquad$\Box$

Now we prove that the  value function satisfy, in the viscosity sense, the
terminal condition.
\begin{lemma}\label{lemTer}
The value function is viscosity solutions to  (\ref{visc2}).\\
\end{lemma}
$Proof.$
First, we prove the subsolution property.

Let us prove that $V^*(T,x)= g(x)$. For this and to begin with, we are going to show that:
\begin{equation}
\begin{array}{ll}
\min\Big[V^*(T,x)-g(x),V^*(T,x)-\sup\limits_{\xi\in U}[V^*(T,x+\xi)-c(T,x,\xi)]\Big]= 0.
\end{array}
\end{equation}
By definition we know that,
$$V^*(T,x)=\varlimsup\limits_{(t',x')\to(T,x),t'<T}V(t',x')\geq\varlimsup\limits_{(t',x')\to(T,x),t'<T}V^{n}(t',x')\qquad\text{for any}\; n\geq 0,$$
therefore
\begin{equation}\label{vgeqG}
V^*(T,x)\geq V^{n}(T,x)=g(x),
\end{equation}
since $V^{n}$ is continuous and at $t = T$ it is equal to $g(x)$. On the other hand we have,
$$V(t,x)\geq \sup\limits_{\xi\in U}[V(t,x+\xi)-c(t,x,\xi)],\qquad\forall (t,x),$$
then
\begin{equation}
V^*(T,x)\geq \sup\limits_{\xi\in U}[V^*(T,x+\xi)-c(T,x,\xi)],
\end{equation}
wich with (\ref{vgeqG}) imply that
\begin{equation}\label{mingeq0}
\min\{V^*(T,x)-g(x),V^*(T,x)-\sup\limits_{\xi\in U}[V^*(T,x+\xi)-c(T,x,\xi)]\}\geq0.
\end{equation}
Let us now show that the left-hand side of (\ref{mingeq0}) cannot be positive.
Let us suppose that for some $x_0$, there is $\epsilon >0$ such that:
\begin{equation}\label{mingeq1}
\min\{V^*(T,x_0)-g(x_0),V^*(T,x_0)-\sup\limits_{\xi\in U}[V^*(T,x_0+\xi)-c(T,x_0,\xi)]\}=2\epsilon.
\end{equation}
Let $(t_k,x_k)_{k\geq 1}$ be a sequence in $[0,T]\times\mathbb{R}^k$ satisfying
$$(t_k,x_k)\to(T,x_0)\qquad\text{and}\qquad V(t_k,x_k)\to V^*(T,x_0)\qquad\text{as}\;k\to\infty.$$
Since $V^*$ is bounded and  $usc$ and taking into account $V^{n}\nearrow V$, we can find
a sequence $(\rho^n)_{n\geq 0}$ of functions of $\mathcal{C}^{1,2}([0,T]\times\mathbb{R}^n)$ such that $\rho^n\to V^*$ and, on
some neighborhood $B_n$ of $(T,x_0)$ we have:
\begin{equation}\label{rho_epsi}
\min\{\rho^n(t,x)-g(x),\rho^n(t,x)-\sup\limits_{\xi\in U}[V^*(t,x+\xi)-c(t,x,\xi)]\}\geq\epsilon,\qquad\forall(t,x)\in B_n.
\end{equation}
After possibly passing to a subsequence of $(t_k,x_k)_{k\geq1}$ we can assume that (\ref{rho_epsi}) holds on $B_n^k:=[t_k,T]\times B(x_k,\delta_k^n)$ for some sufficiently small $\delta^n_k\in(0,1)$ such that $B_n^k\subset B_n$. Now since $V^*$ is  bounded then there exists $\eta> 0$ such that $|V^*| \leq \eta$ on $B_n$. We can then assume that $\rho^n \geq -2\eta$ on $B_n$. Next let us define $\tilde{\rho}^n_k$ by:
$$\tilde{\rho}^n_k(t,x):=\rho^n(t,x)+\frac{4\eta|x-x_k|^2}{(\delta^n_k)^2}+\sqrt{T-t}.$$
Note that $\tilde{\rho}^n_k\geq\rho^n$ and
\begin{equation}\label{vrho}
(V^*-\tilde{\rho}^n_k)(t,x)\leq -\eta\quad\text{for}\;(t,x)\in[t_k,T]\times\partial B(x_k,\delta^n_k).
\end{equation}
As $\partial_t(\sqrt{T-t})\to-\infty$ as $t\to T$, we can choose $t_k$ large enough in front of $\delta^n_k$ and the derivatives of $\rho^n$ to ensure that
\begin{equation}\label{Lrho}
-\mathcal{L}\tilde{\rho}^n_k(t,x)\geq 0 \quad\text{on}\; B_n^k.
\end{equation}
Next let us consider the following stopping time $\theta_n^k:=\inf\{s\geq t_k, (s,X_s^{t_k,x_k})\in {B_n^k}^c\}\wedge T$ where ${B_n^k}^c$ is the complement of $B_n^k$, and $m^k:=\inf\{s\geq t_k, V(s,X_s^{t_k,x_k})=\sup\limits_{\xi\in U}[V(s,X_s^{t_k,x_k}+\xi)-c(s,X_s^{t_k,x_k},\xi)]\}\wedge T.$ Applying now It\^o's formula to the process $(\tilde{\rho}^n_k(s,X_s))_s$ stopped at time $\theta_n^k\wedge m^k$ and taking into account (\ref{rho_epsi}), (\ref{vrho}) and
(\ref{Lrho}) to obtain:
\begin{eqnarray*}
	\tilde{\rho}_k^n(t_k,x_k)
	&\geq&\mathbb{E}\big[\big\{\big(V^*(\theta^n_k,X^{t_k,x_k}_{\theta_n^k})+\eta)\mathbf{1}_{[\theta_n^k< T]}+(\epsilon+g(X^{t_k,x_k}_{T}))\mathbf{1}_{[\theta_n^k= T]}\big\}\mathbf{1}_{[\theta_n^k\leq m^k]}\\ & & +\big\{\epsilon+\sup\limits_{\xi\in U}[V^*(m^k,X_{m^k}^{t_k,x_k}+\xi)-c(m^k,X_{m^k}^{t_k,x_k},\xi)] \big\}\mathbf{1}_{[m^k<\theta_n^k ]}\big] \\ &\geq&
	\mathbb{E}\big[V(\theta_n^k\wedge m^k,X_{\theta_n^k\wedge m^k}^{t_k,x_k})\big]+\eta\wedge\epsilon \\ &=& \mathbb{E}\big[V(t_k,x_k)-\int_{t_k}^{\theta_n^k\wedge m^k}f(s,X_s^{t_k,x_k})ds\big]+\eta\wedge\epsilon.
\end{eqnarray*}
By assumption (H2), we deduce that
$$\lim\limits_{k\to\infty}\mathbb{E}\bigg[\int_{t_k}^{\theta_n^k\wedge m^k}f(s,X_s^{t_k,x_k})ds\bigg]=0.$$
Therefore taking the limit in the previous inequalities yields:
\begin{eqnarray*}
\lim\limits_{k\to\infty}\tilde{\rho}^n_k(t_k,x_k)&=&\lim\limits_{k\to\infty}\rho^n(t_k,x_k)+\sqrt{T-t_k}=\rho^n(T,x_0)\\
&\geq& \lim\limits_{k\to\infty}V(t_k,x_k)+\eta\wedge\epsilon=V^*(T,x_0)+\eta\wedge\epsilon.
\end{eqnarray*}
But this is a contradiction since $\rho^n \to V^*$ pointwisely as $n\to\infty$. Thus for any $x \in \mathbb{R}^k$  we have:
\begin{equation}\label{mingeq2}
\min\{V^*(T,x)-g(x),V^*(T,x)-\sup\limits_{\xi\in U}[V^*(T,x+\xi)-c(T,x,\xi)]\}=0.
\end{equation}
Now, suppose that  $ V^*(T,x)> g(x)$, then from the previous equality there exists $\xi_1\in U$ such that:
\begin{equation}\label{vi=vj}
V^*(T,x)= V^*(T,x+\xi_1)-c(T,x,\xi_1).
\end{equation}
But once more we have $V^*(T,x+\xi_1)>g(x+\xi_1)$. Otherwise, i.e. if $V^*(T,x+\xi_1)=g(x+\xi_1)$ we would have from (\ref{vi=vj}):
\begin{equation}
g(x)<V^*(T,x)= V^*(T,x+\xi_1)-c(T,x,\xi_1)=g(x+\xi_1)-c(T,x,\xi_1),
\end{equation}
which is contradictory with (H4). Therefore
$$V^*(T,x+\xi_1)-\sup\limits_{\xi\in U}[V^*(t,x+\xi_1+\xi)-c(t,x,\xi)]= 0.$$
 Then there exists  $\xi_2\in U$ such that:
\begin{eqnarray*}
	V^*(T,x+\xi_1) &=& V^*(T,x+\xi_1+\xi_2)-c(T,x,\xi_2) \quad\text{and then} \\
	V^*(T,x) &=& V^*(T,x+\xi_1+\xi_2)-c(T,x,\xi_1)-c(T,x,\xi_2).
\end{eqnarray*}
Proceeding in similar fashion, after finitely many steps, boundedness of $V^*$ will be contradicted. Thus we have:
\begin{equation}
\forall x\in\mathbb{R}^k, \quad V^*(T,x)=g(x).\qquad\Box
\end{equation}

\begin{theorem}
	The value function $V,$ defined by (\ref{value function}), is a viscosity solution of the HJB (\ref{eq:HJB})  (with terminal condition $V(T,x) = g(x)$).\\
\end{theorem}
$Proof.$ Thanks to Lemma \ref{lemTer} we have that $V$ satisfies, in the viscosity sense, the terminal condition. As a consequence, we have only to address (\ref{visc1}). First, we prove that $V^*$ is a subsolution to (\ref{visc1}).
Note that since $V^n\nearrow V$ and $V^{n}$ is continuous then we have
\begin{equation}\label{v*}
V^*(t,x)=\lim\limits_{n\to\infty}{\sup}^*V^{n}(t,x)=\varlimsup\limits_{n\to\infty,t'\to t,x'\to x}V^{n}(t',x').
\end{equation}
Since $V^{n}(t,x)\geq\mathcal{M}V^{n-1}(t,x)$, which implies in taking the limit:
\begin{equation}\label{v-geq-h}
\forall (t,x)\in[0,T[\times\mathbb{R}^k \qquad
V^*(t,x)\geq\mathcal{M}V^*(t,x),
\end{equation}
Let us now $(t,x)\in[0,T[\times\mathbb{R}^n$ be such that:
\begin{equation}\label{VMV}
V^*(t,x)-\mathcal{M} V^{*}(t,x)>0.
\end{equation}
Let $(p,q,X)\in \bar{J}^+V^*(t,x)$. By (\ref{v*}) and Lemma 6.1 in \cite{[CIL]}, there exist sequences $$n_j\to\infty, \qquad(p_j,q_j,X_j)\in J^+V^{n_j}(t_j,x_j)$$
such that
$$\lim\limits_{j\to\infty}(t_j,x_j,V^{n_j}(t_j,x_j),p_j,q_j,X_j)=(t,x,V^*(t,x),p,q,X).$$
From the viscosity subsolution property for $V^{n_j}$ at $(t_j,x_j)$, for any $j \geq 0$, we have
\begin{equation}
\min\big\{-p_j-\langle b(t_j,x_j),q_j\rangle-\displaystyle\frac{1}{2}Tr\big[(\sigma\sigma^*)(t_j,x_j)X_j\big] -f(t_j,x_j),V^{n_j}(t_j,x_j)-\mathcal{M} V^{n_j-1}(t_j,x_j)\big\}\leq0,
\end{equation}
Next by (\ref{VMV}), there exists $j_0 \geq 0$, such that if $j \geq j_0$ we have
\begin{equation}
V^{n_j}(t_j,x_j)-\mathcal{M} V^{n_j}(t_j,x_j)>0,
\end{equation}
Therefore for any  $j \geq j_0$
\begin{equation}
-p_j-\langle b(t_j,x_j),q_j\rangle-\displaystyle\frac{1}{2}Tr\big[(\sigma\sigma^*)(t_j,x_j)X_j\big] -f(t_j,x_j)\leq0,
\end{equation}
which implies that
\begin{equation}
-p-\langle b(t,x),q\rangle-\displaystyle\frac{1}{2}Tr\big[(\sigma\sigma^*)(t,x)X\big] -f(t,x)\leq0,
\end{equation}
which is the subsolution property. The supersolution property is proved analogously.$\Box$

Now we give an equivalent of  Hamilton-Jacobi-Bellman equation
(\ref{eq:HJB}).  We consider the new function
$\Gamma$ given by the classical change of variable $\Gamma(t,x) =
\exp(t)V(t, x)$, for any $(t,x)\in[0,T]\times\mathbb{R}^n$.\\ A second property is given by the:
\begin{proposition}
	$V$ is a viscosity solution of (\ref{eq:HJB}) if and only if
	$\Gamma$ is a viscosity solution to the  Hamilton-Jacobi-Bellman equation  in $[0,T)\times \mathbb{R}^n$,
	\begin{equation}
	\begin{array}{ll}\label{eq:HJB1}
 min\Big[\Gamma(t,x)-\cfrac{\partial \Gamma}{\partial t}(t,x)-\mathcal{L}\Gamma(t,x)-\exp(t)f(t,x)\\,\Gamma(t,x)-\widetilde{\mathcal{M}} \Gamma(t,x)\Big]= 0,
	\end{array}
	\end{equation}
	where
	 $$ \widetilde{\mathcal{M}}\Gamma(t,x)=\sup\limits_{\xi\in U}[\Gamma(t,x+\xi)-\exp(t)c(t,x,\xi)],$$
	The terminal
	condition for $\Gamma$ is: $\Gamma(T,x)=\exp(T)g(x)$ in
	$\mathbb{R}^n.$
\end{proposition}
\section{Uniqueness of the solution of Hamilton-Jacobi-Bellman equation}
\no

In this section we deal with the issue of uniqueness of the solution of system (\ref{eq:HJB}) and to do so, we first give the following two  lemmas which is a classical one in viscosity literature (inspired by \cite{[Ik]}):
\begin{lemma}
	\label{con} (convexity of $\mathcal{M}$).
	Let $U,V:[0,T]\times\mathbb{R}^n\rightarrow \mathbb{R}$ and  $\lambda\in[0,1]$. Then
	\begin{eqnarray}
	\mathcal{M}(\lambda U+(1-\lambda)V)&\leq&\lambda \mathcal{M}U+(1-\lambda )\mathcal{M}V\qquad\text{on}\; [0,T]\times\mathbb{R}^n.
	\label{aa}
	\end{eqnarray}
\end{lemma}
	$Proof.$ This follows immediately since
\begin{eqnarray*}	
	\mathcal{M}(\lambda U+(1-\lambda)V)(t,x)&=&\sup\limits_{\xi\in U}[\lambda U(t,x+\xi)+(1-\lambda)V(t,x+\xi)-(\lambda+1-\lambda)c(t,x,\xi)]\\ &\leq& \lambda\mathcal{M}U(t,x)+(1-\lambda)\mathcal{M}V(t,x)
\end{eqnarray*}
	for all $(t,x)\in[0,T]\times\mathbb{R}^n$.\qquad$\Box$
\begin{lemma}
	\label{lemuni}
	Let $V$ a supersolution of $(\ref{eq:HJB1})$.
	Then $V_m:=\big(1-\frac{1}{m}\big)V+\frac{1}{m}\psi$ is a supersolution of
	\begin{equation}
	\label{ttt}
	\begin{array}{ll}
	min\Big[U(t,x)-\cfrac{\partial U}{\partial t}(t,x)-\mathcal{L}U(t,x)-\exp(t)f(t,x)\qquad\\\qquad\qquad\qquad\qquad\qquad\qquad , U(t,x)-\widetilde{\mathcal{M}} U(t,x)\Big]-\displaystyle\frac{k}{m}=0,
	\end{array}
	\end{equation}
	where $ \psi=(\exp(T)\|f\|_\infty +k)$.
\end{lemma}
\\$Proof.$
We have for $(t,x)\in [0,T)\times\mathbb{R}^n$
\begin{equation}
\label{souslemuni}
\begin{array}{ll}
min\Big[\psi(t,x)-\cfrac{\partial \psi}{\partial t}(t,x)-\mathcal{L}\psi(t,x)-\exp(t)f(t,x)\\, \psi(t,x)-\widetilde{\mathcal{M}}\psi(t,x)\Big]\geq min\{\psi-\exp(T)\|f\|_\infty, \exp(t)k\}=k.
\end{array}
\end{equation}
Let $\phi_m\in C^{1,2}([0,T)\times\mathbb{R}^n)$ and  $(\bar{t},\bar{x})\in[0,T)\times\mathbb{R}^n$
such that $\phi_m(\bar{t},\bar{x})=V_m(\bar{t},\bar{x})$, $\phi_m\leq V_m.$ Choose $\phi=\big(\phi_m-\frac{1}{m}\psi\big)\big(\frac{m}{m-1}\big),$ then $\phi(\bar{t},\bar{x})=V(\bar{t},\bar{x})$ and $\phi\leq V$.
As $V$ is a supersolution of $(\ref{eq:HJB1})$, then
\begin{equation*}
\begin{array}{ll}
min\Big[V(\bar{t},\bar{x})-\cfrac{\partial \phi}{\partial t}(\bar{t},\bar{x})-\mathcal{L}\phi(\bar{t},\bar{x})-\exp(\bar{t})f(\bar{t},\bar{x}),V(\bar{t},\bar{x})-\widetilde{\mathcal{M}} V(\bar{t},\bar{x})\Big]\geq 0.
\end{array}
\end{equation*}
Which implies that
\begin{equation*}
\begin{array}{ll}
V_m(\bar{t},\bar{x})-\displaystyle\frac{1}{m}\psi-\displaystyle\frac{\partial}{\partial t}( \phi_{m}-\frac{1}{m}\psi)(\bar{t},\bar{x})-\langle b(\bar{t},\bar{x}),\nabla _x (\phi_{m}-\frac{1}{m}\psi)\rangle  \\+\displaystyle\cfrac{1}{2}tr[\sigma\sigma^*(\bar{t},\bar{x})\nabla^2_x(\phi_{m}-\displaystyle\frac{1}{m}\psi)]-\displaystyle\frac{m-1}{m}f(\bar{t},\bar{x})\geq 0,
\end{array}
\end{equation*}
and then
\begin{equation}\label{uniequi}
\begin{array}{ll}
V_{m}(\bar{t},\bar{x})-\displaystyle\frac{\partial\phi_{m}}{\partial t}(\bar{t},\bar{x})-\langle b(\bar{t},\bar{x}),\nabla _x (\phi_{m})\rangle+\cfrac{1}{2}tr[\sigma\sigma^*(\bar{t},\bar{x})\nabla^2_x(\phi_{m})]-f(\bar{t},\bar{x})\\ -\displaystyle\frac{1}{m}\big(\psi-\cfrac{\partial \psi}{\partial t}-\mathcal{L}\psi-\exp(t)f(\bar{t},\bar{x}) \big) \geq 0.
\end{array}
\end{equation}
Furthermore, by $(\ref{souslemuni})$ we obtain
\begin{equation}
\begin{array}{ll}
V_{m}(\bar{t},\bar{x})-\displaystyle\frac{\partial\phi_{m}}{\partial t} (\bar{t},\bar{x})-\langle b(\bar{t},\bar{x}),\nabla _x (\phi_{m})\rangle+\cfrac{1}{2}tr[\sigma\sigma^*(\bar{t},\bar{x})\nabla^2_x(\phi_{m})]-f(\bar{t},\bar{x})-\displaystyle\frac{k}{m} \geq 0 .
\end{array}
\end{equation}
On the other hand, by the convexity of $\widetilde{\mathcal{M}}$(Lemma \ref{con}), we have
\begin{eqnarray*}
	V_m(\bar{t},\bar{x})-\widetilde{\mathcal{M}}V_m(\bar{t},\bar{x})&\geq& V_m(\bar{t},\bar{x})-\big(1-\frac{1}{m}\big)\widetilde{\mathcal{M}}V(\bar{t},\bar{x})-\frac{1}{m}\widetilde{\mathcal{M}}\psi(\bar{t},\bar{x}) \\ &\geq& V_m(\bar{t},\bar{x})-\big(1-\frac{1}{m}\big)V(\bar{t},\bar{x})-\frac{1}{m}\widetilde{\mathcal{M}}\psi(\bar{t},\bar{x})
	\\ &=&\frac{1}{m}(\psi(\bar{t},\bar{x})-\widetilde{\mathcal{M}}\psi(\bar{t},\bar{x}))> \displaystyle\frac{k}{m}.\qquad\qquad \Box
\end{eqnarray*}

\begin{theorem}\label{comparaison}
	Let $U$ (resp. $V$), be a  bounded u.s.c. viscosity subsolutions (resp.
	bounded l.s.c. viscosity supersolutions) to (\ref{eq:HJB1}). Then, $U\leq V$ in $[0,T]\times\mathbb{R}^n$.
\end{theorem}
\\$Proof.$ We will show by contradiction.
Assume that $\sup_{[0,T]\times\mathbb{R}^n}(U(t,x)-V(t,x))=2M>0$. Furthermore, there exists $m_0$ such that
\begin{equation*}
\sup_{[0,T]\times\mathbb{R}^n}(U(t,x)-V_{m_0}(t,x))>M,
\end{equation*}
where
$V_{m_0}$ is as defined in Lemma $\ref{lemuni}$.

\textbf{Step 1.}

Let $h$ be  a smooth function  such that
\begin{equation}h(x)=\left\{
\begin{array}{c}
0\qquad\qquad\qquad\qquad\qquad\quad\, \text{for}\;|x|\leq 1
\\ \\
\|U\|_{\infty}+\|V_m\|_{\infty}+1\qquad\quad\text{for}\;|x|\geq 2.
\end{array}
\right.
\end{equation}
Then,  for a small $\epsilon, \beta,\gamma > 0$, and $m\geq m_0$ we define:
\begin{equation}
\begin{array}{ll}
\label{phi}
\Phi_{\epsilon,\gamma}(t,x,y)=U(t,x)-V_m(t,y)-\displaystyle\frac{1}{2\epsilon}|x-y|^{2}-\frac{\beta}{T-t}-h_\gamma(x),
\end{array}
\end{equation}
where $h_\gamma(x)=h(\gamma x)$, then we have the following properties:
\begin{itemize}
	\item $h_\gamma(x)=\|U\|_{\infty}+\|V_m\|_{\infty}+1$ when $|x|\geq 2/\gamma$, which ensures that the supremum of $\Phi_{\epsilon,\gamma}$ is achieved therefore is a maximum.
	\item $\nabla h_\gamma(x), D^2h_\gamma(x)\to 0$ as $\gamma\to 0$ uniformly on $\mathbb{R}^n$, which allows to control the differential terms of the $h_\gamma$.
\end{itemize}
Now we consider any maximum points $(t_\epsilon,x_\epsilon, y_\epsilon)$ of the function $\Phi_{\epsilon,\gamma}$.  For $\epsilon$, $\beta$ and $\gamma$ small enough, we have
\begin{eqnarray*}
0<\frac{M}{2}&\leq& \sup_{[0,T]\times\mathbb{R}^n\times\mathbb{R}^n} \Phi_{\epsilon,\gamma}(t,x,y)\leq U(t_\epsilon,x_\epsilon)-V_m(t_\epsilon,y_\epsilon).
\end{eqnarray*}
Next from the inequality
\begin{equation}
\Phi_{\epsilon,\gamma}(t_\epsilon,x_\epsilon,x_\epsilon)+\Phi_{\epsilon,\gamma}(t_\epsilon,y_\epsilon,y_\epsilon)\leq 2\Phi_{\epsilon,\gamma}(t_\epsilon,x_\epsilon,y_\epsilon),
\end{equation}
we get
\begin{eqnarray*}
	\frac{|x_\epsilon-y_\epsilon|^2}{2\epsilon}&\leq& U(t_\epsilon,x_\epsilon)-U(t_\epsilon,y_\epsilon)+V_m(t_\epsilon,x_\epsilon)-V_m(t_\epsilon,y_\epsilon)\quad\\ &=&[U(t_\epsilon,x_\epsilon)-V_m(t_\epsilon,y_\epsilon)]-[U(t_\epsilon,y_\epsilon)-V_m(t_\epsilon,x_\epsilon)].
\end{eqnarray*}
It follows that
\begin{equation}\label{unieps}
|x_\epsilon-y_\epsilon|^2\to 0\qquad \text{as}\;\epsilon \to 0.
\end{equation}
We now want to show that
\begin{equation}\label{uniconvep}
\frac{|x_\epsilon-y_\epsilon|^2}{2\epsilon}\to 0\qquad\text{as}\; \epsilon\to 0.
\end{equation}
Define
$$D:=\sup_{[0,T]\times\mathbb{R}^n}\{U(t,x)-V_m(t,x)
-h_\gamma(x)-\frac{\beta}{T-t}\}.$$
Observe that
\begin{eqnarray*}
	D &\leq& \Phi_{\epsilon,\gamma}(t_\epsilon,x_\epsilon,y_\epsilon)\\&\leq& U(t_\epsilon,x_\epsilon)-V_m(t_\epsilon,y_\epsilon)
	-\frac{\beta}{T-t_\epsilon} -h_\gamma(x_\epsilon)\\ &:=&\phi(\epsilon).
\end{eqnarray*}
By the definition of $\Phi_{\epsilon,\gamma}$, we get (\ref{uniconvep}) if we show that
\begin{equation}\label{uniDphi}
\phi(\epsilon)\to D\qquad \text{as}\;\epsilon\to 0.
\end{equation}
We proceed by contradiction. Assume that (\ref{uniDphi}) does not hold,
there exist $(\epsilon_k,t_{\epsilon_k},x_{\epsilon_k},y_{\epsilon_k})\to (0,\bar{t},\bar{x},\bar{y})$ such that $\lim_{k\to\infty}\phi(\epsilon_k)>D.$
But (\ref{unieps}) implies $\bar{x}=\bar{y}.$ Thus the upper semicontinuity of $U-V_m$ gives
\begin{eqnarray*}
	\lim_{k\to\infty}\phi(\epsilon_k)&\leq& U(\bar{t},\bar{x})-V_m(\bar{t},\bar{x})
	-\displaystyle\frac{\beta}{T-\bar{t}}-h_\gamma(\bar{x})\\ &\leq& D,
\end{eqnarray*}
which gives a contradiction.

\textbf{Step 2. } We now claim that:
\begin{equation}
\label{visco-comp111}
U(t_\epsilon,x_\epsilon)-\sup\limits_{\xi\in U}[U(t_\epsilon,x_\epsilon+\xi)-\exp(t_\epsilon)c(t_\epsilon,x_\epsilon,\xi)]> 0.
\end{equation}
Suppose that
$$U(t_\epsilon,x_\epsilon)-\sup\limits_{\xi\in U}[U(t_\epsilon,x_\epsilon+\xi)-\exp(t_\epsilon)c(t_\epsilon,x_\epsilon,\xi)]\leq 0,$$
then there exist $\xi_1\in U$  such that:
$$U(t_\epsilon,x_\epsilon)-[U(t_\epsilon,x_\epsilon+\xi_1)-\exp(t_\epsilon)c(t_\epsilon,x_\epsilon,\xi_1)]\leq 0.$$
By Lemma \ref{lemuni} we have
$$ V_m(t_\epsilon,y_\epsilon)-V_m(t_\epsilon,y_\epsilon+\xi_1)+\exp(t_\epsilon)c(t_\epsilon,y_\epsilon,\xi_1)\geq \frac{k}{m}.$$
It follows that:
$$\begin{array}{ll}U(t_\epsilon,x_\epsilon)- V_m(t_\epsilon,y_\epsilon) -[U(t_\epsilon,x_\epsilon+\xi_1)
-V_m(t_\epsilon,y_\epsilon+\xi_1)] \\ \leq
-\displaystyle\frac{k}{m}+\exp(t_\epsilon)c(t_\epsilon,x_\epsilon,\xi_1)-\exp(t_\epsilon)c(t_\epsilon,y_\epsilon,\xi_1)\\ \leq
-\displaystyle\frac{k}{m}+\exp(t_\epsilon)\omega_c(|x_\epsilon-y_\epsilon|)
,\end{array}$$
where $\omega_c$ is the modulus of continuity of $c$.Then we have
\begin{eqnarray*}
\Phi_{\epsilon,\gamma}(t_\epsilon,x_\epsilon,y_\epsilon)-\Phi_{\epsilon,\gamma}(t_\epsilon,x_\epsilon+\xi_1,y_\epsilon+\xi_1)&\leq&
h_\gamma(x_\epsilon+\xi_1)-h_\gamma(x_\epsilon)-\displaystyle\frac{k}{m}+\exp(t_\epsilon)\omega_c(|x_\epsilon-y_\epsilon|).	
\end{eqnarray*}
By using the mean value theorem for $h_\gamma$ and by choosing $\epsilon$, $\gamma$ and $m$ appropriately we get
$$	\Phi_{\epsilon,\gamma}(t_\epsilon,x_\epsilon,y_\epsilon)<\Phi_{\epsilon,\gamma}(t_\epsilon,x_\epsilon+\xi,y_\epsilon+\xi).$$
This is contradiction to the fact that $(t_\epsilon, x_\epsilon, y_\epsilon)$ is the supremum point of $\Phi_{\epsilon,\gamma}$. Then the claim (\ref{visco-comp111}) holds.

\textbf{Step 3 } To complete the proof it remains to show
contradiction. Let us denote
\begin{equation}
\varphi_{\epsilon,\gamma}(t,x,y)=\displaystyle\frac{1}{2\epsilon}|x-y|^{2}
+\frac{\beta}{T-t}+h_\gamma(x).
\end{equation}
Then we have: \be \left\{
\begin{array}{lllllll}\label{derive}
	D_{t}\varphi_{\epsilon,\gamma}(t,x,y)=\displaystyle\frac{-\beta}{(T-t)^2}\\
	\nabla_x\varphi_{\epsilon,\gamma}(t,x,y)= \displaystyle\frac{x-y}{\epsilon}+\nabla h_\gamma(x) \\
	\nabla_y\varphi_{\epsilon,\gamma}(t,x,y)= \displaystyle\frac{y-x}{\epsilon}\\
	\\
	B(t,x,s,y)=D^2_{x,y}\varphi_{\epsilon,\gamma}(t,x,y)=\displaystyle\frac{1}{\epsilon}\begin{pmatrix}
		I & -I \\
		-I & I
	\end{pmatrix}+\begin{pmatrix}
		D^2h_\gamma(x) & 0 \\
		0 & 0
	\end{pmatrix}.\\
\end{array}
\right. \ee
Let $c,d\in \mathbb{R}$ such that
$$c + d = \displaystyle\frac{-\beta}{(T-t_\epsilon)^2}.$$
Then applying the result by Crandall et al. (Theorem 8.3, \cite{[CIL]}) to the function
$$U(t,x)-V_m(s,y)-\varphi_{\epsilon,\gamma}(t,x,y)$$
at the point $(t_\epsilon,x_\epsilon,y_\epsilon)$, for any $\epsilon_1>0$, we can find  $X,Y\in S_n,$ such that:
\be \left\{
\begin{array}{lllllll}\label{derive1}
	\big(c,\displaystyle\frac{x_\epsilon-y_\epsilon}{\epsilon}+\nabla h_\gamma(x),X\big)\in J^{2,+}(U(t_\epsilon,x_\epsilon)), \\
	\big(-d,\displaystyle\frac{x_\epsilon-y_\epsilon}{\epsilon},Y\big)\in J^{2,-}(V_m(t_\epsilon,y_\epsilon)),\\
	\\
	-\big(\displaystyle\frac{1}{\epsilon_1}+\|B(t_\epsilon,x_\epsilon,y_\epsilon)\|\big)I\leq \begin{pmatrix}
		X & 0 \\
		0 & -Y
	\end{pmatrix}\leq B(t_\epsilon,x_\epsilon,y_\epsilon)+\epsilon_1B(t_\epsilon,x_\epsilon,y_\epsilon)^2.
\end{array}
\right. \ee
Then by definition of
viscosity solution, we get:
\begin{equation}\begin{array}{lll}\label{vis_sub1}-c+U(t_\epsilon,x_\epsilon)
-\langle\displaystyle\frac{1}{\epsilon}(x_\epsilon-y_\epsilon) +\nabla h_\gamma(x_\epsilon),\\\qquad\qquad
b(t_\epsilon,x_\epsilon)\rangle-\displaystyle\frac{1}{2}tr[\sigma(t_\epsilon,x_\epsilon)\sigma^*(t_\epsilon,x_\epsilon)X]-\exp(t_\epsilon)f(t_\epsilon,x_\epsilon)\leq0,
\end{array}\end{equation}
and
\begin{equation}\begin{array}{l}\label{vis_sub11}d+V_m(t_\epsilon,y_\epsilon)
-\langle\displaystyle\frac{1}{\epsilon}(x_\epsilon-y_\epsilon),\\\qquad\qquad\qquad\qquad
b(t_\epsilon,y_\epsilon)\rangle-\displaystyle\frac{1}{2}tr[\sigma(t_\epsilon,y_\epsilon)\sigma^*(t_\epsilon,y_\epsilon)Y]-\exp(t_\epsilon)f(t_\epsilon,y_\epsilon)\geq \frac{k}{m},
\end{array}\end{equation}
which implies that:
\begin{equation}
\begin{array}{llllll}
\label{viscder} &-c-d+U(t_\epsilon,x_\epsilon)-V_m(t_\epsilon,y_\epsilon)\\ \\& \leq
[\langle\displaystyle\frac{1}{\epsilon}(x_\epsilon-y_\epsilon),
b(t_\epsilon,x_\epsilon)-b(t_\epsilon,y_\epsilon)\rangle\\ \\ &+\langle\nabla h(x_\epsilon), b(t_\epsilon,x_\epsilon)\rangle +\displaystyle\frac{1}{2}tr[\sigma(t_\epsilon,x_\epsilon)\sigma^*(t_\epsilon,x_\epsilon)X-\sigma(t_\epsilon,y_\epsilon)\sigma^*(t_\epsilon,y_\epsilon)Y]\\ \\
&
+\exp(t_\epsilon)f(t_\epsilon,x_\epsilon)-\exp(t_\epsilon)f(t_\epsilon,y_\epsilon)-\displaystyle\frac{k}{m}].
\end{array}
\end{equation}
As
$$ B=B(t_\epsilon,x_\epsilon,y_\epsilon)=\displaystyle\frac{1}{\epsilon}\begin{pmatrix}
I & -I \\
-I & I
\end{pmatrix}+\begin{pmatrix}
D^2h_\gamma(x_\epsilon) & 0 \\
0 & 0
\end{pmatrix}$$
It follows that:
\begin{eqnarray*}
B+\epsilon_1B^2 &\leq& \frac{\epsilon+\epsilon_1}{\epsilon^2}\begin{pmatrix}
I & -I \\
-I & I
\end{pmatrix}+\|D^2h_\gamma(x_\epsilon)\|\begin{pmatrix}
I & 0 \\
0 & 0
\end{pmatrix} + \epsilon_1 \|D^2h_\gamma(x_\epsilon)\|^2\begin{pmatrix}
I & 0 \\
0 & 0
\end{pmatrix},
\end{eqnarray*}
Choosing now $\epsilon_1=\epsilon,$ yields the relation
\begin{eqnarray*}
\label{equaB}
B+\epsilon_1B^2 &\leq& \frac{2}{\epsilon}\begin{pmatrix}
I & -I \\
-I & I
\end{pmatrix}+\|D^2h_\gamma(x_\epsilon)\|\begin{pmatrix}
I & 0 \\
0 & 0
\end{pmatrix}+ \epsilon \|D^2h_\gamma(x_\epsilon)\|^2\begin{pmatrix}
I & 0 \\
0 & 0
\end{pmatrix},
\end{eqnarray*}
Now, from (\textbf{H1}), (\ref{derive1}) and (\ref{equaB}) we get:
\begin{eqnarray*}
	\displaystyle\frac{1}{2}tr[\sigma(t_\epsilon,x_\epsilon)\sigma^*(t_\epsilon,x_\epsilon)X-\sigma(t_\epsilon,y_\epsilon)\sigma^*(t_\epsilon,y_\epsilon)Y]&\leq& \frac{C}{\epsilon}|x_\epsilon-y_\epsilon|^2+C\|D^2h_\gamma(x_\epsilon)\|^2.
\end{eqnarray*}
Next
$$
\langle\frac{1}{\epsilon}(x_\epsilon-y_\epsilon),b(t_\epsilon,x_\epsilon)-b(t_\epsilon,y_\epsilon)\rangle
\leq \frac{C}{\epsilon}|x_\epsilon - y_\epsilon|^{2}.$$ And finally,
\begin{eqnarray*}
	\langle\nabla h_\gamma(x_\epsilon), b(t_\epsilon,x_\epsilon)\rangle\leq
	C\|\nabla h_\gamma(x_\epsilon)\|.
\end{eqnarray*}
So that by plugging into (\ref{viscder}) we obtain:
\begin{equation}
\begin{array}{llllll}
\label{viscder11}
\displaystyle\frac{M}{2}\leq U(t_\epsilon,x_\epsilon)-V_m(t_\epsilon,y_\epsilon)\\ \\
\,\quad\leq \displaystyle\frac{-\beta}{(T-t_\epsilon)^2}+\frac{C}{\epsilon}|x_\epsilon-y_\epsilon|^2+C\|D^2h_\gamma(x_\epsilon)\|^2+C\|\nabla h_\gamma(x_\epsilon)\|
\\ \\  \qquad +\exp(t_\epsilon)f(t_\epsilon,x_\epsilon)-\exp(t_\epsilon)f(t_\epsilon,y_\epsilon)-\displaystyle\frac{k}{m}.
\end{array}
\end{equation}
By sending $\gamma
\rightarrow0$, $\epsilon \rightarrow0$, and $\beta \rightarrow0$, taking into account of the continuity of $f$,
we obtain $\frac{M}{2} \leq 0$, which is a contradiction. Now sending $m\to\infty$, we get the required comparison between $U$ and $V$. The proof of Theorem \ref{comparaison} is now complete. \qquad$\Box$
\begin{corollary}
	Assume that Assumptions (\textbf{H1}), (\textbf{H2}), (\textbf{H3}) and (\textbf{H4}) are fulfilled. Then
	the value function $V$ is a continuous solution of (\ref{eq:HJB}), unique among all bounded solutions.
\end{corollary}

\end{document}